\documentclass [12pt,oneside,mathscr]{amsart}

\usepackage {amscd,amsmath,amssymb,euscript}

\newtheorem{thm}{Theorem}[section]

\newtheorem{lemma}[thm]{Lemma}
\newtheorem{prop}[thm]{Proposition}

\theoremstyle{definition}
\newtheorem{defn}[thm]{Definition}
\newtheorem{example}[thm]{Example}
\newtheorem{examples}[thm]{Examples}

\theoremstyle{remark}
\newtheorem*{pf}{Proof}

\newtheorem*{rmk}{Remark}

\newcommand{\tensor}{\otimes}
\newcommand{\cl}{\operatorname{c}}

\newcommand{\Mod}{\operatorname{Mod}}
\newcommand{\GMod}{\operatorname{G-Mod}}
\newcommand{\SpMod}{\operatorname{Sp-Mod}}
\newcommand{\dual}{\vee}
\renewcommand{\P}{\mathcal P}

\newcommand{\A}{\mathcal A}

\newcommand{\Hom}{\operatorname{Hom}}

\newcommand{\lra}{\longrightarrow}
\newcommand{\lRa}[1]{\stackrel{#1}{\lra}}
\newcommand{\lla}{\longleftarrow}
\newcommand{\lLa}[1]{\stackrel{#1}{\lla}}
\newcommand{\bda}{\bigg{\downarrow}}
\newcommand{\bua}{\bigg{\uparrow}}
\newcommand{\isom}{\cong}
\newcommand{\Z}{{\mathbb Z}}

\newcommand{\QQ}{{\mathcal Q}}

\newcommand{\FF}{{\mathcal F}}

\renewcommand{\H}{\operatorname{H}}
\newcommand{\C}{\mathbb C}

\newcommand{\R}{\mathbf R}

\newcommand{\OO}{{\mathscr O}}
\newcommand{\Ltensor}{\overset{\mathbf L}{\tensor}}

\newcommand{\Ext}{\operatorname{Ext}}

\newcommand{\im}{\operatorname{im}}
\renewcommand{\ker}{\operatorname{ker}}

\newcommand{\D}{\operatorname{D}}

\newcommand{\Xt}{\widetilde{X}}
\newcommand{\Yt}{\widetilde{Y}}
\newcommand{\Pt}{\widetilde{\P}}
\newcommand{\Phit}{\widetilde{\Phi}}

\newcommand{\Psit}{\widetilde{\Psi}}
\newcommand{\Et}{\widetilde{E}}
\newcommand{\Ft}{\widetilde{F}}

\newcommand{\yt}{\tilde{y}}
\newcommand{\xt}{\tilde{x}}

\newcommand{\ft}{\tilde{f}}
\newcommand{\I}{{\mathscr I}}

\newcommand{\eu}{\operatorname{\chi}}

\newcommand{\bSpec}{\operatorname{\mathbf {Spec}}}

\begin{document}
\normalsize
\title[]{Fourier-Mukai transforms for quotient varieties}
\author[]{Tom Bridgeland \and Antony Maciocia}

\date{\today}

\begin{abstract}
We study Fourier-Mukai transforms for smooth projective varieties whose
canonical bundles have finite order. Our results lead to new
transforms for Enriques and bielliptic surfaces.
\end{abstract}

\maketitle

\section{Introduction}
Fourier-Mukai transforms are now well-established as a useful tool for
computing moduli spaces of sheaves on smooth projective varieties [3],
[9]. More recently there has been further interest in these transforms because
of their connection with homological mirror symmetry [8].

In this paper we study Fourier-Mukai transforms for
smooth complex projective varieties whose canonical bundles have finite order, and
relate them to
equivariant transforms on certain finite covering spaces.
Applying our results to the case of Enriques and bielliptic
surfaces, we obtain new examples of transforms for
complex surfaces. These results will be used in [5], where we find all
pairs of minimal complex surfaces with equivalent derived categories.

\smallskip

A Fourier-Mukai (FM) transform is an
exact equivalence
$$\Phi:\D(Y)\lra \D(X)$$
between  the bounded derived categories of coherent sheaves
on two smooth projective varieties $X$ and $Y$.
Due to a result of D. Orlov [14], it is known that for any such equivalence
there is an object $\P$ of $\D(Y\times X)$ and an isomorphism of functors
$$\Phi(-)\isom\R\pi_{X,*}(\P\Ltensor\pi_Y^*(-)),$$
where $Y\lLa{\,\pi_Y}Y\times X\lRa{\pi_X}
X$ are the projection maps.

\smallskip

Let $X$ be a smooth complex projective variety, and suppose that $\omega_X$
has finite order $n$ say. Then there is a finite unbranched cover of $X$ by a smooth
projective variety $\Xt$ with trivial canonical bundle. Moreover, $X$
is the quotient of $\Xt$ by an action of the cyclic group of order
$n$. We call the quotient morphism $$p_X:\Xt\to X$$ the \emph{canonical cover} of $X$.

\smallskip

Suppose $Y$ is another smooth projective variety, and
$$\Phi:\D(Y)\lra\D(X)$$ 
is a FM transform. We show that $\omega_Y$
also has order $n$, and that if $p_Y:\Yt\to Y$ is the canonical cover of
$Y$, there is a $\Z_n$-equivariant FM transform
$$\Phit:\D(\Yt)\lra\D(\Xt)$$ such that the following two squares of functors
commute
$$
\begin{array}{ccc}
\D(\Yt) &\lRa{\Phit} & \D(\Xt) \\
\scriptstyle{p_Y^*}\bua\bda\scriptstyle{p_{Y,*}}
&&\scriptstyle{p_X^*}\bua\bda\scriptstyle{p_{X,*}} \\
\D(Y) &\lRa{\Phi} &\D(X).
\end{array}
$$

Conversely, if $\Yt$ is a smooth projective variety with trivial
canonical bundle, and a $\Z_n$-action with smooth quotient $Y$, and
$$\Phit:\D(\Yt)\lra \D(\Xt)$$
is a $\Z_n$-equivariant FM
transform, then the quotient map $p_Y:\Yt\to Y$ is a canonical cover,
and there is a FM transform
$\Phi$ such that the diagram above commutes.

\subsection*{Notation}

All varieties will be over the complex number field $\C$. Given a projective variety $X$, the category of coherent
$\OO_X$-modules will be denoted $\Mod(X)$. The bounded derived
category of coherent sheaves on $X$ is denoted $\D(X)$. Its objects
are bounded complexes of $\OO_X$-modules with coherent cohomology
sheaves. We refer to [6] for details on derived categories. Note that,
as is usual,
the translation functor on $\D(X)$ is written $[1]$, so that the
symbol $E[m]$ means the object $E$ of $\D(X)$ shifted to the left by
$m$ places.

By a sheaf on $X$ we mean an
object of $\Mod(X)$, and a point of $X$ always means a closed (or
geometric) point. The structure sheaf of such a point $x\in X$ will be
denoted $\OO_x$. The canonical sheaf of a smooth projective variety
$X$ is denoted $\omega_X$.

\section{Canonical covers}

\subsection{}
If $X$ is a smooth projective variety whose canonical bundle has
finite order, one expects a degree $n$ cover of $X$ corresponding to the
element $[\cl_1(\omega_X)]\in\pi_1(X)$. This is the canonical cover of
$X$ referred to in the introduction.

\begin{prop}
\label{bottle}
Let $X$ be a smooth projective variety whose canonical bundle has
finite order $n$. Then there is a smooth projective variety $\Xt$ with
trivial canonical bundle, and
an {\'e}tale cover $p:\Xt\to X$ of degree $n$, such that
\begin{equation} 
\label{lissen}
p_*\OO_{\Xt}\isom \bigoplus_{i=0}^{n-1}\omega_X^i.
\end{equation}
Furthermore, $\Xt$ is uniquely defined up to isomorphism, and there is
a free action of the cyclic group $G=\Z_n$ on $\Xt$ such that
$p:\Xt\to X=\Xt/G$ is the quotient morphism.
\end{prop}

\begin{pf}
By the results of [1], \S I.17, there exists a smooth projective
variety $\Xt$ and a degree $n$ unbranched cover satisfying
(\ref{lissen}). Furthermore $\omega_{\Xt}=p^*\omega_X=\OO_{\Xt}$. By
[7], Ex. II.5.17, $\Xt$ is isomorphic to
$\bSpec(\A)$, where $\A$ is the $\OO_X$-algebra
$$\bigoplus_{i=0}^{n-1}\omega_X^i,$$
which proves uniqueness. The action of $G$ is generated by the
automorphism $\tensor\omega_X$ of $\A$, and clearly $X=\Xt/G$.
\qed
\end{pf} 

\begin{defn}
Let $X$ be a smooth projective variety whose canonical bundle has
finite order $n$. By the \emph{canonical cover} of $X$ we shall mean
the unique smooth projective variety $\Xt$ of Proposition \ref{bottle},
together with the quotient morphism $p_X:\Xt\to X$.
\end{defn}

\begin{examples}

(a) An Enriques surface is a smooth surface $X$ with $\H^1(X,\OO_X)=0$ whose canonical
bundle has order $2$. The canonical cover of such a surface is a K3
surface $\Xt$, and $X$ is the quotient of $\Xt$ by the group generated
by a fixed-point-free automorphism of order $2$. See [1], Ch. VIII.

(b) A bielliptic surface is a smooth surface $X$ with $\H^1(X,\OO_X)=\C^2$ whose
canonical bundle has finite order $n>1$. The possible values of $n$
are $2,3,4$ and $6$. The canonical cover of such a surface is an
Abelian
surface $\Xt$, and $X$ is the quotient of $\Xt$ by a free action of a cyclic group of automorphisms of order
$n$. See [1], \S V.5.
\end{examples}

\subsection{}

Let $X$ be a smooth projective variety whose canonical bundle has
finite order $n$, and let $$p:\Xt\to X$$
 be the canonical cover. Thus $X$ is the quotient of $\Xt$ by a free action of
$G=\Z_n$. Let $g$ be a generator of $G$.

Let $\GMod(\Xt)$ be the category of $G$-equivariant sheaves on
$\Xt$. Since $G$ is
cyclic, a sheaf $E$ on $\Xt$ is equivariant if and only if 
$g^*(E)\isom E$.

Similarly, we let $\SpMod(X)$ denote the category of coherent
$p_*(\OO_{\Xt})$-modules on $X$. A sheaf $E$ on $X$ is in this
category if and only if $E\tensor\omega_X\isom E$. Following [13] we
call these sheaves
\emph{special}.

\begin{lemma}
The functors
$$p^*:\Mod(X)\lra \GMod(\Xt),$$
and
$$p_*:\Mod(\Xt)\lra \SpMod(X),$$
are equivalences of categories.
\end{lemma}

\begin{pf}
This is standard. For the first part see [12], \S 7. The second
part follows from [7], Ex. II.5.17 (e). See also [13], Prop. 1.2.
\qed
\end{pf}

We need to generalise this result to include complexes of sheaves.

\begin{prop}
\label{messy}
(a) Let $\Et$ be an object of $\D(\Xt)$. Then there is an object $E$ of
$\D(X)$ such that $p^*E\isom\Et$ if and only if  there
is an isomorphism $g^*\Et\isom \Et$.

(b)
Let $E$ be
an object of $\D(X)$. Then there is an object $\Et$ of $\D(\Xt)$ such
that $p_*\Et\isom E$ if and only if there is an isomorphism
$E\tensor \omega_X\isom E$.
\end{prop}

\begin{pf}
We shall prove (a); (b) is entirely analogous. One
implication is easy, 
so let us assume that there is
an isomorphism $$s:\Et\lra g^*\Et,$$ and find an object $E$ of $\D(X)$
such that $$p^*E\isom
\Et.$$
We use induction on the number $r$ of non-zero cohomology sheaves of
$\Et$. Shifting $E$ if neccesary, let us assume that $H^i(E)=0$ unless
$-r<i\leq 0$.

The sheaf $\H^0(\Et)$ is $g^*$-invariant, so by the lemma, is isomorphic to
$p^*M$ for some sheaf $M$ on $X$. There is a canonical morphism $\Et\to \H^0(\Et)$,
and hence a triangle
$$\Et\lra \H^0(\Et)\lRa{\ft} \Ft\lra \Et[1],$$
in $\D(\Xt)$, where $\Ft$ has $r-1$ non-zero
cohomology objects. Applying $g^*$ we obtain an isomorphic triangle,
because there is a commutative diagram
$$\begin{CD}
\Et         @>>>         \H^0(\Et) \\
@VsVV                    @V\H^0(s)VV \\
g^*\Et    @>>>         \H^0(g^*\Et).
\end{CD}$$

It follows that $g^*\Ft\isom \Ft$, and so, by induction,
$\Ft\isom p^*F$ for some object $F$ of $\D(X)$.

The lemma below then implies
that $\ft=p^*(f)$ for some morphism $f:M\to F$ of $\D(X)$. Thus there
is an object $E$ of $\D(X)$ and a triangle
$$E\lra M\lRa{f} F\lra E[1].$$
Applying $p^*$ one sees that $p^*E\isom\Et$.
\qed
\end{pf}

\begin{lemma}
Let $M$ be an $\OO_X$-module, and let $F$ be an object of $\D(X)$. Let
$$\ft:p^*M\to p^*F$$ be a morphism of $\D(\Xt)$ such that
$g^*(\ft)=\ft$. Then $\ft=p^*(f)$ for some morphism $f:M\to F$ of
$\D(X)$.
\end{lemma}

\begin{pf}
Replace $F$ by an injective resolution
$$
\begin{CD}
\cdots @>>> I^{-1} @>{d^{-1}}>> I^0 @>{d^0}>> I^1 @>>> \cdots,
\end{CD}
$$
as in [6], Lemma I.4.6.
Then $\ft$ is represented by a morphism
of $\OO_{\Xt}$-modules
$$s:p^*M\to p^*I^0.$$

If $V$ is a finite-dimensional vector space on which $G$ acts, define operators $A$ and $B$ by
$$A=\sum_{i=0}^{n-1} (g^i)^*,\qquad B=1-g^*.$$
Then since $AB=BA=0$, it is easy to check that $\ker A=\im B$.

Take $V$ to be the image of the map
$$p^*(d^{-1})_*:\Hom_{\Xt}(p^*M,p^*I^{-1})\lra\Hom_{\Xt}(p^*M,p^*I^0).$$
The fact that $\ft$ is $G$-invariant means that $B s$ is an element
of $V$. Since $A(B s)=0$, there is an element $k$ of $V$ with
$B k=B s$. Now $$t=s-k\in\Hom_{\Xt}(p^*M,p^*I^0)$$ also represents
$\ft$, and since $B t=0$, is equal to $p^*(u)$ for some
$u\in\Hom_X(M,I^0)$. The result follows.
\qed
\end{pf}
 

\section{Fourier-Mukai transforms}

Let $X$ and $Y$ be smooth projective varieties, and let $\P$ be an
object of $\D(Y\times X)$. Define a functor
$$\Phi^{\P}_{Y\to X}:\D(Y)\lra \D(X)$$
by the formula
$$\Phi^{\P}_{Y\to X}(-)\isom\R\pi_{X,*}(\P\Ltensor\pi_Y^*(-)),$$
where $Y\lLa{\,\pi_Y}Y\times X\lRa{\pi_X}
X$ are the projection maps. Functors of this form will be called
\emph{integral functors}. It is easily checked [10], Prop. 1.3, that
the composite of two integral functors is again an
integral functor.

A Fourier-Mukai (FM) transform relating $X$ and $Y$ is an exact equivalence
of categories
$$\Phi:\D(Y)\lra\D(X).$$
Here exact means commuting with the translation functors and taking
triangles to triangles.
It was proved by Orlov [14], Thm. 2.2, that for any such equivalence there is an object $\P$
of $\D(Y\times X)$, unique up to isomorphism, such that $\Phi$ is
isomorphic to the functor $\Phi^{\P}_{Y\to X}$. We call $\P$ the
\emph{kernel} of the transform $\Phi$.

We shall need the following facts concerning FM
transforms.

\subsection{}
\label{serre}
Recall the definition of a Serre functor on a triangulated category,
[2], pp. 5-6. If $X$ is a smooth projective variety, the functor
$$S_X(-)=(\omega_X\tensor -)[\dim X],$$
is a Serre functor on $\D(X)$. When a Serre functor exists it is
unique up to isomorphism, so any
FM transform relating smooth projective varieties $X$ and  $Y$ must commute
with the functors $S_X$ and $S_Y$. It follows from this that the canonical
bundles of $X$ and $Y$ have the same order.

\subsection{}
If $\P$ is the kernel of a FM transform relating
$X$ and $Y$, then there is an isomorphism
\begin{equation}
\label{glass}
\P\tensor\pi_X^*\omega_X\isom\P\tensor\pi_Y^*\omega_Y.
\end{equation}
Indeed, up to shifts, these objects are the kernels of the left
and right adjoint
functors of $\Phi$ respectively (see e.g. [4], Lemma 4.5), which, since $\Phi$ is an equivalence, must both
be isomorphic to the quasi-inverse of $\Phi$.

\subsection{}
Suppose one has a FM transform $\Phi:\D(Y)\lra\D(X)$ such that for each
$y\in Y$, there is a point $f(y)\in X$ with
$$\Phi(\OO_y)\isom\OO_{f(y)}.$$
I claim that $f$ defines a morphism $Y\to X$, and for some line
bundle $L$ on $Y$, there is an isomorphism of functors
$$\Phi(-)\isom f_*(L\tensor -).$$
To see this note that by [4], Lemma 4.3, the kernel $\P$ of $\Phi$
is a sheaf on $Y\times X$, flat over $Y$, such that for each $y\in Y$,
$\P_y\isom\OO_{f(y)}$. But if $\Delta\subset X\times X$ denotes the
diagonal, then the sheaf $\OO_\Delta$ is a universal sheaf
parameterising structure sheaves of points of $X$. It follows that $f$ is a morphism of varieties,
and $$\P\isom (f\times 1_X)^*(\OO_{\Delta})\tensor\pi_Y^*(L)$$ for some
line bundle $L$ on $Y$. The claim follows.

\subsection{}
Many examples of FM transforms for surfaces are constructed using
the following theorem, which is a simple consequence of the results of
[4]. See [5] for a proof. 

\begin{thm}
\label{freya}
Let $X$ be a smooth projective surface with a fixed polarisation and
let $Y$ be a 2-dimensional, complete,
smooth, fine moduli space of stable, special sheaves on $X$. Then
there is a universal sheaf $\P$ on $Y\times X$ and the resulting
functor $\Phi=\Phi^{\P}_{Y\to X}$ is a FM transform.
\qed
\end{thm}

We shall need the following well-known observation. Suppose we are in the
situation of the theorem, and suppose that $E$ is a stable sheaf on $X$ with the
same Chern character as the sheaves $\P_y$. Then I claim that $E$
must be isomorphic to one of the $\P_y$. If not, for each $y\in Y$, we must have
$$\Hom_{X}(E,\P_y)=0,\qquad\Hom_{X}(\P_y,E)=0.$$
Since $\P_y$ is special, Serre duality implies that
$$\Ext^2_{X}(E,\P_y)=0,$$
and since $E$ has the same Chern character as $\P_y$, and $\Phi$ is an
equivalence
$$\eu(E,\P_y)=\eu(\P_y,\P_y)=\eu(\Phi(\OO_y),\Phi(\OO_y))=\eu(\OO_y,\OO_y)=0,$$
so this is enough to show that $\Hom^i_{X}(E,\P_y)=0$,
for all $i$. 
This is impossible, by [4], Example 2.2, because if $\Psi$ is the quasi-inverse of $\Phi$,
$$\Hom^i_{X}(E,\P_y)=\Hom_{Y}^i(\Psi(E),\OO_y).$$

\subsection{}
We give a couple of well-known examples of FM transforms, which will
be useful later.

\begin{example}
\label{refl}
The first example of an FM transform for a K3 surface was the
reflection functor of [11], although Mukai never explicitly mentions
the fact that it is an equivalence of derived
categories.

To construct it, take a K3 surface $X$ and let $\P$ be the ideal sheaf
$\I_{\Delta}$ 
of the diagonal in $X\times X$. For any $x\in X$,
$\P_x\isom\I_x$. By Theorem \ref{freya}, $\Phi^{\P}_{X\to X}$ is a FM
transform.
\end{example}

\begin{example}
\label{ex2}
Let $(X,\ell)$ be a principally polarised Abelian surface, and let $Y$
be the moduli space of stable sheaves on $X$ of Chern character
$(4,2\ell,1)$. This moduli space is fine,
complete and two-dimensional, so there is a universal sheaf $\P$ on $Y\times X$, and
the resulting functor $\Phi^{\P}_{Y\to X}$ is a FM transform. In fact
$Y$ is isomorphic to $X$. See [9], Prop. 7.1 for details.
\end{example}


\section{Lifts of FM transforms}

In this section we prove our main result, relating FM transforms on
varieties with canonical bundles of finite order, to equivariant
FM transforms on the canonical covers. Throughout we shall suppose that the cyclic group $G=\Z_n$ acts freely on two smooth projective
varieties $\Xt$ and $\Yt$, and denote
the quotient morphisms by $p_X:\Xt\lra X$ and $p_Y:\Yt\lra Y$ respectively.

\begin{defn}
A functor $\Phit:\D(\Yt)\lra \D(\Xt)$ will be called
\emph{equivariant} if there is an automorphism $\mu:G\to G$, and an isomorphism
of functors
$$g^*\circ\Phit\isom \Phit\circ \mu(g)^*,$$
for each $g\in G$.
\end{defn}

\begin{defn}
Given a functor $\Phi:\D(Y)\lra \D(X)$, a \emph{lift} of $\Phi$ is a
functor $\Phit:\D(\Yt)\lra \D(\Xt)$ such that the following
two squares of functors commute up to isomorphism
$$
\begin{array}{ccc}
\D(\Yt) &\lRa{\Phit} & \D(\Xt) \\
\scriptstyle{p_Y^*}\bua\bda\scriptstyle{p_{Y,*}}
&&\scriptstyle{p_X^*}\bua\bda\scriptstyle{p_{X,*}} \\
\D(Y) &\lRa{\Phi} &\D(X),
\end{array}
$$
i.e. such that there are isomorphisms of functors
\begin{equation}
\label{lift}
p_{X,*}\circ\Phit\isom\Phi\circ p_{Y,*},\qquad
p_X^*\circ\Phi\isom\Phit\circ p_Y^*.
\end{equation}
We also say that $\Phit$ \emph{descends} to give the functor $\Phi$.
\end{defn}

We start with a couple of lemmas.

\begin{lemma}
\label{qui}
Let $\Phi:\D(X)\lra \D(X)$ and $\Phit:\D(\Xt)\lra \D(\Xt)$ be integral
functors, such that $\Phit$ lifts $\Phi$.
\begin{enumerate}
\item[(a)] Suppose $\Phi\isom 1_{\D(X)}$. Then $\Phit\isom g_*$ for some
$g\in G$.
\item[(b)] Suppose $\Phit\isom 1_{\D(\Xt)}$. Then $\Phi$ is an
equivalence. If also $p:\Xt\to X$ is the canonical cover, then
$\Phi(-)\isom(\omega_X^{\tensor i}\tensor -)$ for some integer $i$.
\end{enumerate}
\end{lemma}

\begin{pf}
We start with (a). Take a point $\xt\in\Xt$, and put $x=p_X(\xt)$. Then $E=\Phit(\OO_{\xt})$ satisfies
$p_{X,*}(E)=\OO_x$, so $E=\OO_{f(\xt)}$ for some
point $f(\xt)$ in the fibre $p^{-1}(x)$. By (3.3),
$f:\Xt\to\Xt$ is a morphism of varieties, and for some line bundle $L$
on $\Xt$,
$$\Phit(-)\isom f_*(L\tensor-).$$
Since $f(\xt)$ always lies
in the fibre $p^{-1}(x)$, $f=g$ for some $g\in G$. Now the functor $g^{-1}_*\circ\Phit$
also lifts the identity, and takes $p_X^*(\OO_X)=\OO_{\Xt}$ to $L$, so
in fact $L$ must be trivial.

To prove (b), take a point $x\in X$, and a point $\xt\in\Xt$ such that
$p_X(\xt)=x$. Then
$$\Phi(\OO_x)=p_{X,*}(\OO_{\xt})=\OO_x.$$
It follows from (3.3) that $\Phi\isom (L\tensor -)$ for some line bundle $L$
on $X$. We must have $p_X^*L=\OO_{\Xt}$, so if $p_X$ is the canonical
cover, the projection formula gives
$$L\tensor(\bigoplus_{i=0}^{n-1}\omega_X^i)=L\tensor p_{X,*}(p_X^*\OO_X)=p_{X,*}(p_X^*L)=\bigoplus_{i=0}^{n-1}\omega_X^i,$$
hence $L$ is a power of $\omega_X$. 
\qed
\end{pf}

\begin{lemma}
\label{rightback}
Let $\Pt$ and $\P$ be objects of $\D(\Yt\times\Xt)$ and $\D(Y\times X)$
respectively, such that
\begin{equation}
\label{drag}
(p_Y\times 1_X)^*(\P)\isom (1_{\Yt}\times p_X)_*(\Pt).
\end{equation}
Then $\Phit=\Phi^{\Pt}_{\Yt\to\Xt}$ is a lift of
$\Phi=\Phi^{\P}_{Y\to X}$.
\end{lemma}

\begin{pf}
Put
$$f=(1_{\Yt}\times p_X),\qquad h=(p_Y\times 1_X),$$
and consider the commutative diagram
$$
\begin{array}{ccccc}
\Yt &\lLa{\pi_{\Yt}} &\Yt\times\Xt &\lRa{\pi_{\Xt}} &\Xt \\
\parallel &          &\bda\scriptstyle{f} && \bda\scriptstyle{p_X} \\
\Yt &\lLa{j} &\Yt\times X &\lRa{k} & X \\
\scriptstyle{p_Y}\bda && \scriptstyle{h}\bda &&\parallel \\
Y &\lLa{\pi_Y} &Y\times X &\lRa{\pi_X} & X.
\end{array}
$$
Let $E$ be an object of $\D(\Yt)$. By [6], II.5.6, II.5.12, there are natural isomorphisms
\begin{eqnarray*}
p_{X,*}(\Phit(E))& = & p_{X,*}\R\pi_{\Xt,*}(\Pt\Ltensor\pi_{\Yt}^*E) \\
   &\isom & \R k_* (f_*(\Pt\Ltensor f^* j^* E)))\isom\R k_*(f_*\Pt\Ltensor j^*E) \\
   &\isom & \R\pi_{X,*} (h_*(h^*\P\Ltensor j^* E)) \isom\R\pi_{X,*} (\P\Ltensor h_* j^* E) \\
   &\isom & \R\pi_{X,*} (\P\Ltensor \pi_Y^*(p_{Y,*} E))= \Phi(p_{Y,*}(E)).
\end{eqnarray*}
The second isomorphism of (\ref{lift}) can be proved in the same way,
or by taking adjoints.

\qed
\end{pf}

The following theorem is the main result of this paper.

\begin{thm}
\label{bec}
Let $X$ and $Y$ be smooth projective varieties with canonical bundles
of order $n$, and take canonical covers 
$$p_X:\Xt\to X, \qquad p_Y:\Yt\to
Y.$$
Thus $X$ and $Y$ are quotients of $\Xt$ and $\Yt$  by
the cyclic group $G=\Z_n$. Then any FM transform 
\begin{equation}
\label{ttwo}
\Phi:\D(Y)\lra \D(X)
\end{equation}
lifts to give an
equivariant FM transform
\begin{equation}
\label{oone}
\Phit:\D(\Yt)\lra \D(\Xt).
\end{equation}
Conversely, any equivariant FM transform (\ref{oone}) is
the lift of some  FM transform (\ref{ttwo}).
\end{thm}

\begin{pf}
First let $\Phi:\D(Y)\lra \D(X)$ be a FM transform, and let $\P$ be
its kernel. Put 
$$\QQ=(p_Y\times 1_X)^*(\P).$$
It follows from the isomorphism
(\ref{glass}) that $\QQ\tensor \omega_{\Yt\times
X}\isom \QQ$, so by Prop. \ref{messy}, there is an object $\Pt$ of $\D(\Yt\times\Xt)$ satisfying
(\ref{drag}). Define $\Phit=\Phi^{\Pt}_{\Yt\to\Xt}$. Then by Lemma
\ref{rightback}, $\Phit$ is a
lift of $\Phi$.

Let $\Psi$ be a quasi-inverse for $\Phi$. Then $\Psi$ is also an FM
transform and hence lifts to a
functor $\Psit:\D(\Yt)\lra \D(\Xt)$ by the same argument. Now it is easy
to check that $\Psit\circ\Phit$ is a lift of $\Psi\circ\Phi\isom
1_{\D(Y)}$. Hence, by Lemma \ref{qui}, composing $\Psit$ with
$g^*$ for some $g\in G$, we can assume that $\Psit\circ\Phit\isom
1_{\D(\Yt)}$. Similarly, $\Phit\circ\Psit\isom 1_{\D(\Xt)}$, so
$\Phit$ is an equivalence.

Take $g\in G$ and consider the FM transform
$g^*\circ \Phit$. This is also a lift of $\Phi$, so
$\Psit\circ g^*\circ \Phit$ is a lift of $1_{\D(Y)}$. By Lemma \ref{qui} again, there is an
element $\mu(g)\in H$ such that
$$g^*\circ \Phit\isom\Phit\circ(\mu(g))^*.$$
Clearly, the homomorphism $\mu:G\to G$ must be injective, and so by symmetry it is an
isomorphism.

For the converse, let $\Phit:\D(\Yt)\lra\D(\Xt)$ be a FM transform
with kernel $\Phit$. The $G$-equivariance of $\Phit$ is equivalent to the condition
$$(1_{\Yt}\times g)^*(\Pt)\isom (\mu(g)\times 1_{\Xt})^*(\Pt)\quad\forall
g\in G.$$
It follows that $(1_{\Yt}\times p_X)_*(\Pt)$ is $G$-invariant so that
$$(p_Y\times 1_X)^*(\P)\isom (1_{\Yt}\times p_X)_*(\Pt)$$
for some object $\P$ of $\D(Y\times X)$. Hence by Lemma \ref{rightback}, $\Phit$
lifts $\Phi=\Phi^{\P}_{Y\to X}$.

We must show that $\Phi$ is an equivalence of categories. Let $\Psit$
be a quasi-inverse of $\Phit$. Then $\Psit$ is $G$-equivariant and hence
is the lift of some integral functor $\Psi:\D(X)\lra \D(X)$. But then $\Psit\circ\Phit\isom
1_{\D(\Yt)}$ lifts $\Psi\circ\Phi$, so by Lemma \ref{qui},
twisting $\Psi$ by some power of $\omega_X$,
$\Psi\circ\Phi\isom 1_{\D(Y)}$. Similarly $\Phi\circ\Psi\isom
1_{\D(X)}$. 

\qed
\end{pf}

\begin{rmk}
In the situation of the
theorem, it is easy to see using Lemma \ref{qui} that if two FM transforms
$\Phit_1,\Phit_2$ lift a given transform $\Phi$, then $\Phit_2\isom
g^*\circ\Phit_1$ for some $g\in G$.

Similarly, if FM transforms
$\Phi_1,\Phi_2$ both lift to give the same transform $\Phit$, then
$\Phi_2\isom \omega_X^i\tensor\Phi_1$ for some integer $i$.
\end{rmk}

\smallskip

A couple of points remain. Let $X$ be a smooth projective variety $X$ whose canonical bundle
has order $n$, and let $p_X:\Xt\to X$ be the canonical cover. Thus $X$
is the quotient of $\Xt$ by a free action of $G=\Z_n$.

Firstly,
suppose there is a FM transform $\Phi$ relating $X$ to another
variety $Y$. Then by (\ref{serre}), $\omega_Y$ also has order $n$, and taking canonical
covers of $X$ and $Y$ we are in the situation of Theorem \ref{bec}.

Secondly, suppose there is another smooth projective variety $\Yt$
with a free $G$-action, and that there is an equivariant FM transform
$\Phit$ relating $\Xt$ and $\Yt$. Then I claim that the quotient
morphism $p_Y:\Yt\to Y$ is a canonical cover of $Y=\Yt/G$, so we are
again in the situation of Theorem \ref{bec}.

To prove the claim, note that by the argument used in the proof of
Theorem \ref{bec}, the functor $\Phit$ descends to give a FM transform
$\Phi:\D(Y)\lra \D(X)$. By the result of (3.1) $\omega_Y$ has order $n$. Taking a canonical cover $Y'$ of $Y$ we can
lift $\Phi$ to a FM transform $\Phi':\D(Y')\lra\D(\Xt)$. Now
$\Phit^{-1}\circ\Phi'$ is an equivariant FM transform relating $Y'$ and $\Yt$ which
lifts the identity on $\D(Y)$. It follows that $Y'$ and $\Yt$ are
isomorphic as $G$-spaces, hence the claim.


\section{Examples}

Let $\Xt$ be a smooth projective surface with a fixed polarisation and let $\Yt$ be a complete, fine, smooth,
two-dimensional moduli space of stable sheaves on $\Xt$. Then there is
a universal sheaf $\Pt$ on $\Yt\times\Xt$, and by Theorem \ref{freya}, the resulting
functor
$$\Phit=\Phi^{\Pt}_{\Yt\to\Xt}:\D(\Yt)\lra\D(\Xt),$$
is a FM transform.

Assume that $\Xt$ has trivial canonical bundle (so is of either
Abelian or K3 type). As we noted in (3.1), $\Yt$ also has
trivial canonical bundle. Suppose
further that the cyclic group $G=\Z_n$ acts freely on $\Xt$ via
automorphisms. Let $p_X:\Xt\to X$ denote the quotient morphism.

Applying the result of (3.3) it is easy to see that there is an algebraic action of $G$ on the moduli space $\Yt$
such that for each point $\yt\in\Yt$ and each $g\in G$,
\begin{equation}
\label{can}
g^*(\Pt_{\yt})\isom\Pt_{g(\yt)}.
\end{equation}
If the action of $G$ on $\Yt$ is free then we can form the quotient
$Y=\Yt/G$, and Theorem \ref{bec} shows that $\Phit$ descends to give a FM transform
$\Phi:\D(Y)\lra\D(X)$. The following lemma gives a purely numerical
criterion for when this happens.

\begin{lemma}
\label{dump}
The action of $G$ on $\Yt$ defined above is free, if and only if the
highest common factor of the integers
$$\eu(p_X^* F,\Pt_{\yt})=\eu(p_X^* F^{\dual}\tensor \Pt_{\yt}),$$
as $F$ varies through all vector bundles on $X$ is 1.
\end{lemma}

\begin{pf}
Let $E=\Pt_{\yt}$. Note first that by the adjunction $p_X^*\dashv p_{X,*}$
$$\eu(p_X^* F,E)=\eu(F,p_{X,*}E).$$
If the action of $G$ on $\Yt$ is free then as we noted above $\Phit$
descends to
a transform $\Phi:\D(Y)\lra\D(X)$. Then if $y=p_Y(\yt)$,
$\Phi(\OO_y)=p_{X,*}E,$ so if $\Psi$ is the inverse of $\Phi$,
$$\eu(F,p_{X,*}(E))=\eu(\Psi(F),\OO_y).$$
Since $\Psi$ is an equivalence the highest common factor of these
integers is 1.

For the converse let $g$ be a generator of $G$, and suppose that the $G$-action is not free, so that for
some proper factor $m$ of $n$, the element $g^m$ of $G$ fixes $E$. Then the sheaf
$$\bigoplus_{i=0}^{m-1} (g^i)^*(E),$$
is $g^*$-invariant, so by Prop. \ref{messy}, is isomorphic to $p_X^* A$ for some sheaf $A$ on
$X$. 

For any bundle $F$ on $X$,
$$m\eu(p_X^* F,E)=\eu(p_X^* F,p_X^* A)=\eu(p_{X,*}
p_X^*(F),A)=n\eu(F,A),$$
because
$$p_{X,*}
p_X^*(F)\isom F\tensor(\bigoplus_{i=0}^{n-1}\omega_X^i).$$
It follows that $n/m$ divides $\eu(p_X^* F,E)$.
\qed
\end{pf}

\begin{example}
Let $X$ be an Enriques surface. Then there is a K3 surface $\Xt$ with an
automorphism $\sigma$ of order 2 such that $X$ is the quotient of $\Xt$ by the
2-element group generated by $\sigma$. For any point $\xt\in \Xt$ one
has
$$\sigma^*(\I_{\xt})=\I_{\sigma(\xt)},$$
so the reflection functor of Example \ref{refl} descends to give an FM
transform
$$\Phi:\D(X)\lra \D(X),$$
This has the property that for each $x\in X$ one has an exact sequence
$$0\lra\Phi(\OO_x)\lra \OO_X\oplus \omega_X\lra \OO_x\lra 0.$$
It is this transform which was studied in [15], \S 3.7.
\end{example}

\begin{example}
Let $X$ be a bielliptic surface whose fundamental group is cyclic of
order $n$. Then the canonical cover of $X$ is a product of elliptic curves
$\Xt=C_1\times C_2$, and $X$ is the quotient of $\Xt$ by a free action
of $G=\Z_n$.

The original Fourier-Mukai functor of [10] never descends
because the sheaf
$\OO_{\Xt}=\FF(\OO_0)$ is $G$-invariant.

Consider instead the moduli space $\Yt$ of stable sheaves on $\Xt$ of Chern
character $(4,2\ell,1)$, where $\ell=C_1+C_2$ is a principal
polarisation.  By Lemma \ref{dump}, the FM transform of Example \ref{ex2} descends to give an FM
transform
$$\Phi:\D(Y)\lra \D(X),$$
such that $\Phi(\OO_y)$ is a locally free sheaf of rank $4n$ for all
$y\in Y$.
\end{example}

\begin{rmk}

Any bielliptic surface $X$ is the quotient of a product of elliptic
curves $C_1\times C_2$ by some finite Abelian group $G$, but the quotient
map $C_1\times C_2\to X$ is only the canonical cover of $X$ if $G$ is
cyclic. Thus in general, a FM transform $\D(X)\lra \D(X)$ will \emph{not}
lift to a transform $\D(C_1\times C_2)\lra \D(C_1\times C_2)$.
\end{rmk}

In [5] we shall show that if $X$ and $Y$ are Enriques or bielliptic
surfaces, and $\Phi:\D(Y)\lra \D(X)$ is a FM transform, then $X$ and
$Y$ are isomorphic.

\section*{References}

\small

[1] {\it W. Barth, C. Peters, A. Van de Ven,} Compact Complex
Surfaces, Ergebnisse Math. Grenzgeb. (3), vol. 4, Springer-Verlag,
1984.

[2] {\it A. Bondal, D. Orlov,} Reconstruction of a variety from the
derived category and groups of autoequivalences, Preprint alg-geom 9712029.

[3] {\it T. Bridgeland,} Fourier-Mukai transforms for elliptic
surfaces, J. reine angew.
math. {\bf 498} (1998) 115-133, also alg-geom 9705002.

[4] {\it T. Bridgeland,} Equivalences of triangulated categories and
Fourier-Mukai transforms, to appear in
Bull. Lond. Math. Soc. (1998), also alg-geom 9809114.

[5] {\it T. Bridgeland, A. Maciocia,} Complex surfaces with equivalent
derived categories, in preparation.

[6] {\it R. Hartshorne,} Residues and duality, Lect. Notes Math. {\bf
20}, Springer-Verlag, 1966.

[7] {\it R. Hartshorne,} Algebraic Geometry, Grad. Texts Math. {\bf 52},
Springer-Verlag, 1977.

[8] {\it M. Kontsevich,} Homological algebra of mirror symmetry, Preprint
alg-geom 9411018.

[9] {\it A. Maciocia,} Generalized Fourier-Mukai transforms, J. reine
angew. Math. {\bf 480} (1996), 197-211, also alg-geom 9705001.

[10] {\it S. Mukai,} Duality between $D(X)$ and $D(\Hat X)$ with its
application to Picard sheaves, Nagoya Math. J. {\bf 81} (1981),
153-175.

[11] {\it S. Mukai,} On the moduli space of bundles on K3 surfaces I,
in: Vector Bundles on Algebraic Varieties, M.F. Atiyah et al., Oxford
University Press (1987), 341-413.

[12] {\it D. Mumford,} Abelian varieties, Oxford University Press (1970).

[13] {\it D. Naie,} Special rank two vector bundles over Enriques
surfaces, Math. Ann. {\bf 300} (1994) 297-316.

[14] {\it D.O. Orlov,} Equivalences of derived categories and K3
surfaces, J. Math. Sci. (NY), {\bf 84} (1997) 5 1361-1381, also
alg-geom 9606006.

[15] {\it S. Zube,} Exceptional vector bundle on Enriques surfaces,
Mat. Zametki {\bf 61} (1997) 6 825-834 also alg-geom 9410026.
 
\bigskip

Department of Mathematics and Statistics, The University of Edinburgh,
King's Buildings, Mayfield Road, Edinburgh, EH9 3JZ, UK.

\smallskip

email: {\tt tab@maths.ed.ac.uk}
{\tt\ \ \ \ \ \ \ \ \ \ ama@maths.ed.ac.uk}

\end{document}